\documentclass[12pt]{article}
\usepackage{latexsym,color,amsmath,amsthm,amssymb,amscd,amsfonts}

\newtheorem{lemma}{Lemma}[section]

\newtheorem{theorem}[lemma]{Theorem}
\newtheorem{definition}[lemma]{Definition}
\newtheorem{corollary}[lemma]{Corollary}

\newtheorem{conj}[lemma]{Conjecture}

\newtheorem*{remark*}{Remark}

\def\R{{\mathbb R}}

\def\qed{\hfill $\vcenter{\hrule height .3mm
\hbox {\vrule width .3mm height 2.1mm \kern 2mm
\vrule width .3mm height 2.1mm} \hrule height .3mm}$ \bigskip}

\begin{document}
\title{The
$M$-ellipsoid, Symplectic Capacities and Volume }
\author{Shiri Artstein-Avidan, Vitali Milman\thanks{ The first named author was supported by the National
Science Foundation under agreement No.~DMS-0111298. The first and second
named authors were supported in part by a grant from the
US-Israeli BSF.} \, and {Yaron} Ostrover}
\maketitle
\noindent {\it {\bf Abstract:}
In this work we bring together tools and ideology from two
different fields,
Symplectic Geometry and Asymptotic Geometric Analysis, to arrive at some
new results. Our main result is a dimension-independent bound for
the symplectic capacity of a convex body by its volume radius.}

\section{Short Introduction}
In this work we bring together tools and ideology from two different
fields, Symplectic Geometry and Asymptotic Geometric Analysis, to
arrive at some new results. Our main result is a
dimension-independent bound for the symplectic capacity of a convex
body by its volume radius. This type of inequality was first
suggested by C. Viterbo, who conjectured that among all convex
bodies in ${\R}^{2n}$ with a given volume, the Euclidean ball has
maximal symplectic capacity (definition in Section~\ref{SGB} below).
More precisely, Viterbo's conjecture states that the best possible
constant $\gamma_n$ such that for any choice of a symplectic
capacity $c$ and any convex body $K \subset {\R}^{2n}$ we have
\[  {{\frac {c(K)} {c(B^{2n})}} \leq \gamma_n \biggl ({\frac {{\rm
Vol}(K)} {{\rm Vol}(B^{2n})}} \biggr )^{1/n}}, \]
is $\gamma_n = 1$, where $B^{2n}$ is the Euclidean unit ball in
${\mathbb R}^{2n}$. The estimate which Viterbo proved in his work
\cite{V} was $\gamma_n\le {32n}$, and in the case of centrally
symmetric bodies he showed that $\gamma_n\le {2n}$. Hermann showed
in~\cite{H} that for the special class of convex Reinhardt domains
the conjecture holds. The first and third named authors showed
in~\cite{AO} that there exists a universal constant $A_1$ such
that $\gamma_n\le A_1 (\log 2n)^2$, and also presented wide
classes of bodies where the inequality holds without the
logarithmic term. The methods used for that result came from
Asymptotic Geometric Analysis. In this work we use some more
advanced methods of Asymptotic Geometric Analysis to show that the
logarithmic term is not needed at all. That is, there exists a
universal constant $A_0$ for which $\gamma_n \le A_0$ for any
dimension $2n$ and all convex bodies in $\R^{2n}$. Moreover, these
strong bounds are obtained by using linear tools only. While this
fits with the philosophy of Asymptotic Geometric Analysis, this is
less expected from the point of view of Symplectic Geometry where
for strong results one expects to need highly nonlinear objects.

\noindent {\bf Notations:} In this paper the letters
$A_0,A_1,A_2,A_3$ and $C$ are used to denote universal positive
constants which do not depend on the dimension nor on the body
involved.
In what follows we identify ${\mathbb R}^{2n}$ with ${\mathbb C}^n$
by associating to $z = x + iy$, where $x,y \in {\mathbb R}^n$, the
vector $(x_1,y_1,\ldots,x_n,y_n)$, and consider the standard complex
structure given by complex multiplication by $i$, i.e.
$i(x_1,y_1,\ldots,x_n,y_n) = (-y_1,x_1,\ldots,-y_n,x_n)$. We denote
by $\langle \cdot , \cdot \rangle$ the standard Euclidean inner
product on ${\mathbb R}^{2n}$.
We shall 
denote by $e^{i\theta}$ the standard action of $S^1$ on ${\mathbb
C}^n$ which rotates each coordinate by angle $\theta$, i.e., $e^{i
\theta}(z_1,\ldots,z_n)= (e^{i \theta}z_1,\ldots,e^{i
\theta}z_n)$. By $x^{\perp}$ we denote the hyperplane orthogonal
to $x$ with respect to the Euclidean inner product. For two sets
$A, B$ in ${\mathbb R}^{2n}$, we denote their Minkowski sum by
$A+B = \{a+b : a \in A, b \in B \}$. By a convex body we shall
mean a convex bounded set in $\R^{2n}$ with non-empty interior.
Finally, since affine translations in ${\mathbb R}^{2n}$ are
symplectic maps, we shall assume throughout the text that any
convex body $K$ has the origin in its interior.

\noindent {\bf Structure of the paper:} The paper is organized as
follows. In the next section we recall the necessary definitions
from symplectic geometry, describe the history of the problem and
state our main theorem. In Section~\ref{AGB} we describe the main
tool coming from Asymptotic Geometric Analysis, called the
$M$-ellipsoid. In Section~\ref{proofsec} we prove our main result,
and in the last section we show an additional result about convex
bodies, generalizing a result of Rogers and Shephard.

\noindent {\bf Acknowledgments:} The third named author thanks
Leonid Polterovich for his kind support and for helpful advice
regarding the text. The first and the second named authors thank the
Australian National University where part of this work was carried
out.

\section{Symplectic geometry background and the main result.} \label{SGB}

Consider the $2n$-dimensional Euclidean space ${\mathbb R}^{2n}$
with the standard linear coordinates $(x_1,y_1, \ldots,x_n,y_n)$.
One equips this space with the standard symplectic structure
$\omega_{st} = \sum_{j=1}^n dx_j \wedge dy_j$, and with the
standard inner product $g_{st} = \langle \cdot,\cdot \rangle$.
Note that under the identification 
between ${\mathbb R}^{2n}$ with $ {\mathbb C}^n$ these two
structures are the real and the imaginary part of the standard
Hermitian inner product in ${\mathbb C}^n$, and
$\omega(v,iv) = \langle v , v \rangle $.

In~\cite{V}, Viterbo related between the symplectic way of
measuring the size of sets using what is called ``symplectic
capacities", and the classical Riemannian approach, using the
canonical volume. Among other things, he conjectured that in the
class of convex bodies in ${\mathbb R}^{2n}$ with fixed volume,
the Euclidean ball has maximal symplectic capacity. This
isoperimetric inequality was proved in the same paper~\cite{V} up
to a constant $\gamma_n$ which is linear in the dimension (see
Theorem~\ref{cnjV} below).

In this work, we continue the approach taken in~\cite{AO}, where
methods from Asymptotic Geometric Analysis were used to reduce the
order of the above mentioned constant $\gamma_n$. In \cite{AO} it
was improved from order $n$ to order $(\log n)^2$, where $n$ is
the dimension of the space involved. In this note we improve it
further, and prove an upper bound for $\gamma_n$ which is
independent of the dimension. Finding dimension independent
estimates is a frequent goal in Asymptotic Geometric Analysis,
where surprising phenomena such as concentration of measure (see
e.g.~\cite{MS}) imply the existence of order and structures in
high dimension despite the huge complexity it involves. It is
encouraging to see that such phenomena also exist in Symplectic
Geometry, and although this is just a first example, we hope more
will follow. Furthermore, we wish to stress that the tools we use
are purely linear and the reader should not expect any difficult
symplectic analysis.

In order to state our results we continue with the formal
definitions.


\begin{definition} \label{Def-sym-cap}
A
symplectic capacity on $({\mathbb
R}^{2n},\omega_{st})$ associates to each  subset $U \subset
{\mathbb R}^{2n}$ a non-negative number $c(U)$ such that the
following three properties hold:
\begin{enumerate}
\item[(P1)] $c(U) \leq c(V)$ for $U \subseteq V$ (monotonicity)
\item[(P2)] $c \big (\psi(U) \big )= |\alpha| \, c(U)$ for  $\psi
\in {\rm Diff} ( {\mathbb R}^{2n} )$ such that $\psi^*\omega_{st}
= \alpha \, \omega_{st}$ (conformality)
\item[(P3)] $c \big (B^{2n}(r) \big ) = c \big (B^2(r) \times {\mathbb
C}^{n-1} \big ) = \pi r^2$ (nontriviality and normalization),
\end{enumerate}
\end{definition}
\noindent where $B^{2k}(r)$ is the open $2k$-dimensional ball of
radius $r$. Note that the third property disqualifies any
volume-related invariant, while the first two properties imply
that every two sets $U,V \subset {\mathbb R}^{2n}$ will have the
same capacity provided that there exists a symplectomorphism
sending $U$ onto $V$. Recall that a {\it symplectomorphism} of
${\mathbb R}^{2n}$ is a diffeomorphism which preserves the
symplectic structure i.e., $\psi \in {\rm Diff} ( {\mathbb R}^{2n}
)$ such that $\psi^* \omega_{st} = \omega_{st}$. We will denote by
${\rm Symp}({\mathbb R}^{2n}) = {\rm Symp}({\mathbb
R}^{2n},\omega_{st})$ the group of all symplectomorphisms of
$({\mathbb R}^{2n},\omega_{st})$.

A priori, it is not clear that symplectic capacities exist. The
celebrated non-squeezing theorem of Gromov~\cite{G} shows that for
$R
> r$ the ball $B^{2n}(R)$ does not admit a symplectic embedding
into the symplectic cylinder $Z^{2n}(r):= B^2(r) \times {\mathbb
C}^{n-1}$. This theorem led to the following definitions:

\begin{definition} The symplectic radius of a non-empty set
$U \subset {\mathbb R}^{2n}$ is
$$ c_B(U) := \sup \left \{\pi r^2 \, | \,
\ There \ exists \  \psi \in {\rm Symp}({\mathbb R}^{2n}) \ with \
\psi \left (B^{2n}(r) \right ) \subset U \right  \}.$$ The
cylindrical capacity of $U$ is
$$ {c}^Z(U) := \inf \left \{\pi r^2 \, | \,
\ There \ exists \  \psi \in {\rm Symp}({\mathbb R}^{2n}) \ with \
\psi (U) \subset Z^{2n}(r)  \right  \}.$$
\end{definition}

Note that both the symplectic radius and the cylindrical capacity
satisfy the axioms of Definition~\ref{Def-sym-cap} by the
non-squeezing theorem. Moreover, it follows from
Definition~\ref{Def-sym-cap} that for every symplectic capacity
$c$ and every open set $U \subset {\mathbb R}^{2n}$ we have
$c_B(U) \le c(U) \le c^Z(U)$.

The above axiomatic definition of symplectic capacities is
originally due to Ekeland and Hofer~\cite{EH}. Nowadays, a variety
of symplectic capacities can be constructed in different ways. For
several of the detailed discussions on symplectic capacities we
refer the reader
to~\cite{CHLS},~\cite{Ho},~\cite{HZ},~\cite{L},~\cite{Mc}
and~\cite{V1}.

In this work we are interested in an inequality relating the
symplectic capacity of a convex body in ${\R}^{2n}$ and its
volume. Viterbo's conjecture states that
among all convex bodies in ${\R}^{2n}$ with a given volume,
the symplectic capacity is maximal for the Euclidean ball. Note
that by monotonicity this is obviously true for the symplectic
radius $c_B$. More precisely, denote by ${\rm Vol}(K)$ the volume
of $K$ and abbreviate $B^{2n}$ for the open Euclidean unit ball in
${\mathbb R}^{2n}$. Following Viterbo~\cite{V} we state

\begin{conj} \label{The-conjecture} For any symplectic capacity
$c$ and for any convex
body $K \subset {\mathbb R}^{2n}$
$${\frac {c(K)} {c(B^{2n})}} \leq
 \biggl ( {\frac {{\rm Vol}(K)} {{\rm Vol}(B^{2n})} } \biggr
)^{1/n}$$ and equality is achieved only for symplectic images of
the Euclidean ball.
\end{conj}

The first result
in this direction is due to Viterbo~\cite{V}. Using John's
ellipsoid he proved:
\begin{theorem}[Viterbo]\label{cnjV}
For a convex body $K \subset {\mathbb R}^{2n}$ and a symplectic
capacity $c$ we have
$$ {\frac {c(K)} {c(B^{2n})}} \leq
\gamma_n   \biggl ( {\frac {{\rm Vol}(K)} {{\rm Vol}(B^{2n})} }
\biggr )^{1/ n}$$ where $\gamma_n = {2n}$ if $K$ is centrally
symmetric and $\gamma_n = 32n$ for general convex bodies.
\end{theorem}

In \cite{H}, Hermann constructed starshaped domains in ${\mathbb
R}^{2n}$, for $n > 1$, with arbitrarily small volume and fixed
cylindrical capacity. Therefore, in the category of starshaped
domains the above theorem with {\em any} constant $\gamma_n$
independent of the body $K$ must fail. In addition, he proved the
above conjecture for a special class of convex bodies which admit
many symmetries called convex Reinhardt domains (for definitions
see \cite{H}).

In \cite{AO}, the first and third named authors showed
\begin{theorem}
There exists a universal constant $A_1$ such that
for a convex body $K \subset {\mathbb R}^{2n}$ and a symplectic
capacity $c$ we have
\[ {\frac {c(K)} {c(B^{2n})}} \leq
A_1 (\log 2n)^2  \biggl ( {\frac {{\rm Vol}(K)} {{\rm Vol}(B^{2n})} }
\biggr )^{1/ n}\]
(so, $\gamma_n \le A_1 (\log {2n})^2$).
\end{theorem}
They also showed that for many classes of convex bodies, the
logarithmic term is not needed.  Among these classes are all the
$\ell_p^n$-balls for $1\le p\le \infty$, all zonoids (bodies that
can be approximated by Minkowski sums of segments) and other wide
classes of convex bodies, see \cite{AO}.

In this work we eliminate the logarithmic factor from the above
theorem. Before we state our main results we wish to re-emphasize
that, as in \cite{AO}, we work exclusively in the category of linear
symplectic geometry. That is, we restrict ourselves to the concrete
class of linear symplectic transformations. It turns out that even
in this limited category, the tools are powerful enough to obtain a
dimension independent estimate for $\gamma_n$ in Theorem~\ref{cnjV}.
More precisely, let ${\rm Sp}({\mathbb R}^{2n})={\rm Sp}({\mathbb
R}^{2n},\omega_{st})$ denote the group of linear symplectic
transformation of ${\mathbb R}^{2n}$. We consider a more restricted
notion of linearized cylindrical capacity, which is similar to $c^Z$
but where the transformation $\psi$ is taken only in ${\rm
Sp}({\mathbb R}^{2n})$ namely
$$ {c}_{lin}^Z(U) := \inf \left \{\pi r^2 \, | \,
{\rm \ There \ exists } \  \psi \in {\rm Sp}({\mathbb R}^{2n}) \
{\rm with} \ \psi (U) \subset Z^{2n}(r)  \right  \}.$$ Of course,
it is always true that for every symplectic capacity $c$ we have
$c \le c^Z \le c_{lin}^Z$.

Our main result is that for some universal constant $A_0$ one has
that $\gamma_n \le A_0$ for all $n$. This follows from the following
theorem, which we prove in Section \ref{proofsec}.
\begin{theorem}\label{TMT}
There exists a universal constant $A_0$ such that for every even
dimension $2n$ and any convex body $K \subset {\mathbb R}^{2n}$
\[  {\frac {c^Z_{lin}(K)} {c(B^{2n})}}
\leq A_0
 \biggl (
{\frac {{\rm Vol} (K)} {{\rm Vol} (B^{2n})} } \biggr )^{ 1/ n}.\]
\end{theorem}


\section{Asymptotic geometric analysis background: $M$-position} \label{AGB}

In this section we work in $\R^n$ with the Euclidean structure,
without a symplectic or complex structure. We review some well
known theorems from Asymptotic Geometric Analysis which we will
use in later sections. A position of a convex body is equivalent
to a choice of a Euclidean structure, or, in other words, a choice
of some ellipsoid as the Euclidean unit ball. A fundamental object
in Asymptotic Geometric Analysis, which was discovered by the
second named author in relation with the reverse Brunn-Minkowski
inequality, is a special ellipsoid now called the Milman
ellipsoid, abbreviated $M$-ellipsoid. This ellipsoid has several
essentially equivalent definitions, the simplest of which may be
the following:

\begin{definition}
An ellipsoid ${\cal E}_K$ is called an $M$-ellipsoid (with constant $C$)
of $K$ if
${\rm Vol}({\cal E}_K) = {\rm Vol}(K)$ and it satisfies
\[{\rm Vol} (K+{\cal E}_K)^{1/n} \le C{\rm Vol}(K)^{1/n}, \quad {\rm and} \quad
{\rm Vol} (K\cap {\cal E}_K)^{1/n}\ge C^{-1}{\rm Vol}(K)^{1/n}.\]
 \end{definition}

The fact that there exists a universal $C$ such that every convex
body $K$ has an $M$-ellipsoid (with constant $C$) was proved in
\cite{Mil} for a symmetric body $K$. The fact that the body $K$
need not be symmetric, for the existence of an $M$-ellipsoid with
the properties
%
which we use in the proof of our main result, was proved in
\cite{Milm1317} (see Theorem 1.5 there). A complete extension of
all $M$-ellipsoid properties in the non-symmetric case was
performed in \cite{MilPaj}, where it was shown that the right
choice of the origin (translation) in the case of a general convex
body is the barycenter (center of mass) of the body.

This ellipsoid was invented in order to study
the ``reverse Brunn-Minkowski inequality" which is proved
in \cite{Mil}, and we begin by recalling this
inequality, which we will strongly use in the proof of
our main theorem. We then describe some further properties of this ellipsoid.
Recall that the classical Brunn-Minkowski
inequality states that if $A$ and $B$ are non-empty compact
subsets of ${\mathbb R}^n$, then
$$ {\rm Vol} (A+B)^{1/n} \geq {\rm Vol}(A)^{1/n} + {\rm
Vol}(B)^{1/n}.$$ Although at first sight it seems that one cannot
expect any inequality in the reverse direction (imagine, for
example, two very long and thin ellipsoids pointing in orthogonal
directions in $\R^2$), if one allows for an extra choice of
``position'', a reverse inequality is possible.

It was discovered in \cite{Mil} that one can reverse the
Brunn-Minkowski inequality, up to a universal constant factor, as
follows: for every convex body $K$ there exists a linear
transformation $T_K$, which is volume preserving, such that for
any two bodies $K_1$ and $K_2$, the bodies $T_{K_1}K_1$ and
$T_{K_2}K_2$ satisfy an inverse Brunn-Minkowski inequality. It
turned out that the right choice of $T_K$ is such that the
ellipsoid $r(T_K)^{-1}B^n$ (for the right choice of $r$) is an
$M$-ellipsoid of $K$, which we denote as before by  ${\cal E}_K$.
We then say that the body $T_K K$ is in $M$-position (or that it
is an $M$-position of $K$). Thus, a body is in $M$-position if a
multiple of the Euclidean ball $B^n$ is an $M$-ellipsoid for $K$.
We remark that an $M$-ellipsoid of a body is far from being
unique, and a body can have many different such ellipsoids. For a
detailed account about $M$-ellipsoids we refer the reader to
\cite{Milm1317} and \cite{Pi1}, where there are also proofs of the
theorems below. The property of $M$-position which we use in this
paper for the proof of Theorem \ref{TMT} is the following

\begin{theorem} \label{RBM}
There exists a universal constant $C$ such that if ${\widetilde
K_1}, {\widetilde K_2} \subset \R^n$ are two convex bodies in
$M$-position then
\begin{equation}\label{heya}  {\rm Vol}({\widetilde K_1}
+ {\widetilde K_2})^{1/n} \leq C \left ( {\rm Vol}({\widetilde
K_1})^{1/n} + {\rm Vol}({\widetilde K_2})^{1/n} \right ).
\end{equation}
\end{theorem}

In particular this theorem implies that for a convex body $K$
there exists a transformation $T_K$, which depends solely on $K$,
such that
%
for any two convex bodies $K_1$ and $K_2$, denoting ${\widetilde
K_1} = T_{K_1}(K_1), {\widetilde K_2} = T_{K_2}(K_2)$, we have
that (\ref{heya}) is satisfied. The transformation $T_K$ is the
transformation which takes the ellipsoid ${\cal E}_K$ to a
multiple of $B^{n}$. Therefore, it is clear that any composition
of $T_K$ with an orthogonal transformation from the left will also
satisfy this property.

This ellipsoid ${\cal E}_K$ has many more well known intriguing properties.
We recall one of them, which we will use in Section \ref{appli}:
\begin{theorem}\label{withPP}
There exists a universal constant $C$ such that for any convex body
$K$, the ellipsoid ${\cal E}_K$ satisfies the following: for every
convex body $P$ 
one has that
\begin{equation}\label{withP}
C^{-1}{\rm Vol}(P + {\cal E}_K)^{1/n} \le
{\rm Vol}(P + K)^{1/n} \le  C{\rm Vol}(P + {\cal E}_K)^{1/n}. \end{equation}
\end{theorem}

\section{Proof of the Main Result}\label{proofsec}


We return to ${\mathbb R}^{2n}$ equipped with the standard
symplectic structure and the standard Euclidean inner product. We
first present the main ingredient needed for the proof of the main
theorem. With aid of the $M$-position, we show that every convex
body $K$ has a linear symplectic image $K' = SK$ such that the
couple $K'$ and $iK'$ satisfy the inverse Brunn-Minkowski
inequality. For this we need to recall a well known fact about the
relation between a symplectic form and a positive definite
quadratic form.  The following theorem by Williamson~\cite{W}
concerns simultaneous normalization of a symplectic form and an
inner product.

\noindent{\bf Williamson's theorem:} {\it For any positive
definite symmetric matrix $A$ there exists an element $S \in {\rm
Sp}(2n)$ and a diagonal matrix with positive entries $D$  with the
property $iD = Di$ (complex linear), such that $A = S^TDS$.}

An immediate corollary (for a proof see~\cite{AO}) is

\begin{corollary} \label{decomposition of the
m-position} Let $T$ be a volume preserving $2n$-dimensional real
matrix. Then there exists a linear symplectic matrix $S \in {\rm
Sp}({\mathbb R}^{2n})$, an orthogonal transformation $W \in O(2n)$
and a diagonal complex linear matrix $D$ with positive entries
such that
\[ T  = WDS. \]
\end{corollary}

This decomposition, together with Theorem~\ref{RBM}, implies the
following (in the sequel we will only use the special case $\theta
= \pi/2$, i.e., multiplication by $i$)

\begin{theorem} \label{SIK}
Every convex body $K$ in ${\mathbb R}^{2n}$ has a
symplectic image $K'=SK$, where $S \in {\rm Sp}(2n)$,
such that for any $0\le \theta\le 2\pi$
\[  {\rm Vol}(K)^{1/2n} \leq {\rm Vol}(K'+e^{i\theta}K')^{1/2n} \le A_2 {\rm
Vol}(K)^{1/2n},\]
where $A_2$ is a universal constant.
\end{theorem}

\noindent{\bf Proof.} The first inequality holds trivially for any
$K' = SK$ since $K' \subset K' + e^{i \theta}K'$. Next, let $K$ be
a convex body in ${\mathbb R}^{2n}$. Set $K_1= TK$, where $T$ is a
volume-preserving linear transformation which takes the body $K$
to an $M$-position. It follows from Corollary~\ref{decomposition
of the m-position} that $T =WDS$ where $W$ is orthogonal, $S$ is
symplectic, and $D$ is a complex linear transformation. We set $K'
= SK$. The remark after Theorem~\ref{RBM} implies that we can
assume $K_1 = DSK$ where $D$ and $S$ are as above, since an
orthogonal
image 
of a body in $M$-position is also in $M$-position.
Note that the rotated body $e^{i\theta}
 K_1$ is in $M$-position as well, since multiplication
by a complex number of module 1 is a unitary transformation. Next,
it follows from Theorem \ref{RBM} that
\[ {\rm Vol}(K_1 + e^{i\theta}
K_1)^{1/2n} \leq C \left ( {\rm Vol}(K_1)^{1/2n} + {\rm
Vol}(e^{i\theta}K_1)^{1/2n} \right ) = 2C {\rm Vol}(K)^{1/2n} ,\]
where $C> 0$ is a universal constant. Since $D$ is complex linear it
commutes with multiplication by $e^{i\theta}$, and using also the
fact that it is volume preserving we conclude that
\[ {\rm Vol} (K' + e^{i\theta}K')^{1/2n} = {\rm Vol} (K_1 +
e^{i\theta}K_1)^{1/2n} \le 2C {\rm Vol}(K)^{1/2n}.\] The proof is
now complete. $\hfill \square$

In order to complete the proof of the main theorem, we shall need
two more ingredients. The first is the following easy observation
\begin{lemma}\label{ai}
Let $K$ be a symmetric convex body satisfying $K = iK$, and let 
$rB^{2n} \subset K$ be the largest multiple of the Euclidean
ball contained in $K$. Then
\[ c^Z_{lin}(K) \le 2\pi r^2. \] 
\end{lemma}

\noindent{\bf Proof.} 
Since the body $K$ is assumed to be symmetric there are at least
two contact points $x$ and $-x$ which belong to $\partial K$, the
boundary of $K$, and to $rS^{2n-1}$, the boundary of $rB^{2n}$.
Note that the supporting hyperplanes to $K$ at these points must be
$\pm x + x^{\perp}$ since they are also supporting hyperplanes of
$rB^{2n}$ at the tangency points. Thus, the body $K$ lies between the
hyperplanes  $-x + x^{\perp}$ and $x + x^{\perp}$. However, since
$K$ is invariant under multiplication by $i$, the points $\pm ix$
are contact points for $\partial K$ and $rS^{2n-1}$ as well. Thus,
the body $K$ lies also between $-ix + ix^{\perp}$ and $ix +
ix^{\perp}$. Note that the length of the vectors $x$ and $ix$ is
$r$. We conclude that the projection of $K$ onto the plane spanned
by $x$ and $ix$ is contained in a square of edge length $2r$, which
in turn is contained in a disc of radius $\sqrt{2}r$. Therefore $K$
is contained in a cylinder of radius $\sqrt{2}r$ with base spanned
by $x$ and $ix$. Since this cylinder is a unitary image of the
standard symplectic cylinder $Z^{2n}(\sqrt{2}r)$ the lemma follows.
$\hfill \square$

\noindent{\bf Remark:} The factor $2\pi$ above can be replaced by
$4$ if we replace $c^Z_{lin}$ by $c^Z$. For this we need only to
take a small step out of the linear category and use a non-linear
symplectomorphism which is essentially two-dimensional.

The last tool we need is a famous result of Rogers and
Shephard~\cite{RogShe}. This result, which we generalize in some
sense in Section \ref{appli} below, states that for a convex body
$K \subset \R^n$ the volume of the so called ``difference body''
$K-K$ is not much larger than the volume of the original body.
They show that one has
\begin{equation}\label{RS}
{\rm Vol}(K-K) \le 4^{n}{\rm Vol}(K).\end{equation}

We are now in a position to prove our main result:

\noindent{\bf Proof of Theorem \ref{TMT}.} Let $K$ be a convex
body in $\R^{2n}$ and set $K_1 = K - K$. Note that $K_1$ is
symmetric and by (\ref{RS}) we have
${\rm Vol}(K_1) \le 4^{2n}{\rm Vol}(K)$. It follows from Theorem~\ref{SIK}
 that there exists a symplectic map $S \in {\rm Sp}({\mathbb R}^{2n})$ for which ${\rm
Vol}(SK_1+iSK_1) \le A_2^{2n} {\rm Vol}(K_1)$. Denote $K_2 =
SK_1$, $K_3 = K_2 + iK_2$.
Thus ${\rm Vol}(K_2) = {\rm Vol}(K_1)$ and
${\rm Vol}(K_3) \le A_2^{2n} {\rm Vol}(K_2)$.
Let $r>0$ be the largest radius
such that $rB^{2n} \subset K_3$. We thus have
\[ r^{2n} {\rm Vol} (B^{2n}) \le
{\rm Vol} (K_3) \le  A_2^{2n} {\rm Vol}(K_2) =  A_2^{2n} {\rm
Vol}(K_1) \le (4A_2)^{2n} {\rm Vol}(K). \] On the other hand,
since $K_3 = iK_3$, it follows from the monotonicity property of
symplectic capacities and from Lemma~\ref{ai} that
\[ {c^Z_{lin}(K)}  \le  {c^Z_{lin}(K_1)} =
{c^Z_{lin}(K_2)} \le {c^Z_{lin}(K_3)} \le  2\pi r^2. \]

Joining these two together we conclude
\[  \frac{c^Z_{lin}(K)}{c(B^{2n})} \le 2(4A_2)^2
\left( \frac{{\rm Vol}(K)}{ {\rm Vol} (B^{2n})}\right)^{1/n}, \]
and the proof of the theorem is complete.  $\hfill \square$

\section{Generalized Rogers Shephard}\label{appli}

In this section we again work in $\R^n$ equipped only with the
Euclidean structure.
The above type of  reasoning led us to the following simple
generalization of the theorem of Rogers and Shephard (\ref{RS})
above. In this generalization, instead of considering the Minkowski
sum and the Minkowski difference of a body and itself, we consider
the sum and the difference of two different bodies, and show with the
use of $M$-ellipsoid that both have the same volume radius up to a
universal constant. We remark that the constant in (\ref{RS}) is
equal to $2$ (if we put it in the setting of the theorem below)
whereas
the constant in the theorem below, although universal, may be
much worse.

\begin{theorem}\label{GRS}
There exists a universal constant $A_3$ such that
for any two convex bodies $A, B\subset \R^n$ one has
\[ {\rm Vol}(A+B)^{1/n} \le A_3 {\rm Vol}(A-B)^{1/n} . \]
\end{theorem}

\noindent{\bf Proof.} In the case where one of the bodies is centrally
symmetric the statement is trivial.
In the case where both of them are not symmetric, we will use
the property of the $M$-ellipsoid described in Theorem
\ref{withPP} above.
Let ${\cal E}_B$ be the $M$-ellipsoid of $B$, which is of
course centrally symmetric.
We see that
\[ {\rm Vol}(A+B)^{1/n} \le C {\rm Vol}(A+{\cal E}_B)^{1/n}
=  C {\rm Vol}(A-{\cal E}_B)^{1/n} \le
C^{2} {\rm Vol}(A-B)^{1/n}. \]
$\hfill \square$

\end{document}